\documentclass{article}
\usepackage{amsmath,amsopn,amssymb}
\allowdisplaybreaks

\newcommand{\bb}{\begin{bmatrix}}
\newcommand{\eb}{\end{bmatrix}}
\numberwithin{equation}{section}
\newtheorem{Theorem}{Theorem}[section]
\newtheorem{Lemma}{Lemma}[section]

\newtheorem{Remark}{Remark}[section]

\newtheorem{Corollary}{Corollary}[section]

\catcode`@=11
\newbox\@temp
\def \blap#1{\vbox to 0mm{#1\vss}}

\newcommand{\bset}[3][2mm]{#3\llap{%
      \blap{\vskip#1\hbox to 0mm{\hss $#2$\hss}}%
      \setbox\@temp\hbox{\mbox{$#3$}}\hskip0.5\wd\@temp}}
\def\hhline{%
  \noalign{\ifnum0=`}\fi\hrule \@height 3\arrayrulewidth \futurelet
   \@tempa\@xhline}
\catcode`@=12

\title{A note on the $\top$-Stein matrix equation
}

\author{Chun-Yueh Chiang
\thanks{Corresponding Author, Center for General Education,National Formosa University, Huwei 632, Taiwan. {\tt (chiang@nfu.edu.tw)}}
}

\begin{document}
\maketitle

\begin{abstract}
This note is concerned with the linear matrix equation
$X = AX^\top B + C$, where the operator $(\cdot)^\top$ denotes
the transpose ($\top$) of a matrix. The first part of this paper set forth the necessary and sufficient conditions for the unique solvability of the solution $X$.
The second part of this paper aims to provide a comprehensive treatment of the relationship between the theory of the generalized eigenvalue problem and the theory of the linear matrix equation. In the finally part of this paper starts with a briefly review of numerical methods for solving the linear matrix equation. Related to the computed methods, knowledge of the residual is discussed.  An expression related to the backward error of an approximate solution is obtained; it shows that a small backward error implies a small residual. Just like for the discussion of linear matrix equations, perturbation bounds for solving the linear matrix equation are also proposed in this work.

\end{abstract}


\pagestyle{myheadings} \thispagestyle{plain}

\section{Introduction}
Our purpose of this work is to study the so-called $\top$-Stein matrix equation
\begin{equation}\label{DTS}
X = AX^\top B+C,
\end{equation}
where $A,\, B,\, C \in\mathbb{R}^{n\times n}$ are known matrices, and $X\in \mathbb{R}^{n\times n}$ is an unknown matrix to be determined.
Our interest in the $\top$-Stein equation originates from
the study of completely integrable mechanical systems, that is, the
analysis of the $\top$-Sylvester equation
\begin{equation}\label{CTS}
A X + X^\top B = C,
\end{equation}
where $A$, $B$, $C$ are matrices in $\mathbb{R}^{n\times n}$~\cite{Braden98, Ma2003}.
By means of the generalized inverses
or QZ decomposition~\cite{Bojanczyk1992}, the solvability conditions of~\eqref{CTS} are studies in~\cite{Braden98, Ma2003, Chiang2012}.
Suppose that the matrix pencil $A-\lambda B^\top$ is regular, that is,
$ a A + b  B^\top$ is invertible for some scalars $a$ and $b$.
The $\top$-Sylvester equation~\eqref{CTS} can be written as\begin{equation}\label{add}
(aA + b B^\top)X + X^\top(aB+bA^\top) = aC+ b C^\top.
\end{equation}
Pre-multiplying both sides of~\eqref{add} by $(aA + b B^\top)^{-1}$, we have
\begin{equation}
X + UX^\top V = D,
\end{equation}
where $U = (aA + b B^\top)^{-1}$, $V = aB+bA^\top$ and
$ D =  (aA + b B^\top)^{-1}(aC+ b C^\top)$. This is of the form~\eqref{DTS}.
In other words, numerical approaches for solving~\eqref{CTS} can be obtained by
transforming~\eqref{CTS} into the form of~\eqref{DTS}, and then applying numerical methods to~\eqref{DTS} for the solution~\cite{Chiang2012,Su2010, Wang2007}.
With this in mind, in this note we are interested in the study of
$\top$-Stein matrix equation~\eqref{DTS}.


Our major purpose in this work can be divided into three parts. First, we determine necessary and sufficient conditions for the unique solvability of the solution to~\eqref{DTS}.
In doing so, Zhou et al.~\cite{Zhou2011} transform~\eqref{DTS} to
the standard Stein equation
\begin{equation}\label{Stein}
W = AB^\top W A^\top B + AC^\top B +C
\end{equation}
with respect to the unknown matrix $W\in\mathbb{R}^{m\times n}$
and give the following necessary condition
\begin{equation}\label{condition}
\mu \nu  \neq 1,\quad \mbox{for all } \mu, \nu \in \sigma(A^\top B).
\end{equation}
Here, $\sigma(A^\top B)$ be the set of all eigenvalues of $A^\top B$. Zhou shows that if~\eqref{Stein} has a unique solution, then~\eqref{DTS} has a unique solution. However,
 a counterexample is provided in~\cite{Zhou2011} to show that the relation~\eqref{condition} is only a necessary condition for the unique solvability of~\eqref{DTS}.

In~\cite{Chiang2012, Ikramov2012}, the periodic QZ (PQZ) decomposition~\cite{Bojanczyk1992} is applied to consider the necessary and sufficient conditions of the unique solvability of~\eqref{DTS},
conditions are given in~\cite{Ikramov2012} ignore the possibility of the existence of the unique solution, while $1$ is a simple root of
$\sigma(A^\top B)$. This condition is included in our subsequent discussion and the following remark is provided to support
our observation.
\begin{Remark}\label{rem1}
Let $A = -1$ and $B = 1$, that is, $\sigma (AB^\top) = \{-1\}$. It is clear that, the scalar equation $X=-X^\top+C $ has a unique solution $X=\dfrac{C}{2}$. But, condition~\eqref{condition} is not satisfied by choosing $\mu = \nu = -1$.
\end{Remark}
It can also be observed from Remark~\ref{rem1} that even if \eqref{DTS} is uniquely solvable, it does not imply~\eqref{Stein} (namely, $X=X+C-C$) is uniquely solvable.
Conditions in~\cite[(4.6)]{Chiang2012} provided that conditions for the unique solvability of the solution to~\eqref{DTS} via a structured algorithm. In our work, we through a complete analysis
 for square coefficient matrices in terms of the analysis of the spectra of the matrix $A^\top B$, the new approach to the condition of unique solvability of
the $\top$-Stein equation~\eqref{DTS} can be obtained.

Second, we present the invariant subspace method and, more
generally, the deflating subspace method to solve the $\top$-Stein equation. Our
methods are based on the analysis of the eigeninformation for a matrix pencil. We carry out a thorough discussion to address the
various eigeninformation encountered in the subspace methods. These ideas can be
implemented into algorithms easily.

Finally, we take full account of the error analysis of~ Eq.~\eqref{DTS}. Expressions and implications such as the residual, the backward error, and perturbation bounds are derived in this work. Note that for an approximate solution $Y$ of~\eqref{DTS}, the backward error tells us
how much the matrices $A$, $B$ and $C$ must be perturbed. An important point found in Section~\ref{sec:error_ana} is that a small backward error indicates a small value for the residual  $\mathcal{R} = Y-AY^\top B - C$, but reverse is not usually true.

Beginning in Section 2, we formulate the necessary and sufficient conditions for the existence of the solution of~\eqref{DTS} directly by means of the spectrum analysis. In Section 3 we provide an deflating subspace method for computing the solution of Eq.~\eqref{DTS}. Numerical methods for solving Eq.~\eqref{DTS} and the related residual analysis are discussed in Section 4. The associated error analysis of Eq.~\eqref{DTS} is given in Section 5 and concluding remarks are given in Section 6.

%
%
\section{Solvability conditions of the Matrix Equation \eqref{DTS}  } 
In order to formalize our discussion, let the notations $A\otimes
B$ be the Kronecker product of matrices $A$ and $B$, $I_n$ be the $n\times n$ identity matrix, and $\| \cdot \|_F$ denotes the Frobenius norm.

With the
Kronecker product, Eq.~\eqref{DTS} can be written as
the enlarged linear system
\begin{align}\label{KronD}
(I_{n^2}-(B^\top\otimes A)\mathcal{P}) \mbox{vec}(X)=\mbox{vec}(C),
\end{align}
where $\mbox{vec}(X)$ stacks the columns of $X$ into a column vector
and $\mathcal{P}$ is the Kronecker permutation matrix~\cite{Bernstein2009} which maps
$\mbox{vec}(X)$ into $\mbox{vec}(X^\top)$, i.e.,
\begin{equation*}
\mathcal{P}=
\sum\limits_{1\leq i,j\leq n}e_je_i^\top \otimes e_ie_j^\top,
\end{equation*}
where $e_i$ denotes the $i$-th column of the $n\times n$ identity
matrix $I_{n}$. Due to the specific structure of $\mathcal{P}$, it
has been shown in~\cite[Corollary 4.3.10]{Horn1994} that
\begin{equation*}
\mathcal{P}^\top (B^\top\otimes A) \mathcal{P}  = A\otimes B^\top.
\end{equation*}
It then follows that
\begin{equation}\label{kroeig}
 ((B^\top\otimes A) \mathcal{P} )^2= (B^\top\otimes A)\mathcal{P}\mathcal{P}^\top (A\otimes B^\top)   = B^\top A\otimes A B^\top,
\end{equation}
since $\mathcal{P}^2 = I_{n^2}$ and $\mathcal{P} = \mathcal{P}^\top$.
Note that eigenvalues of matrices $A^\top B$ and $A B^\top$ are the same.
By~\eqref{kroeig} and the property of the Kronecker
product~\cite[Theorem 4.8]{zhang1999}, we know that
\begin{equation*}
\sigma(((B^\top\otimes A) \mathcal{P})^2) = \left \{ \lambda_i
\lambda_j | \lambda_i, \lambda_j \in \sigma(A^\top B) = \left\{
\lambda_1,\ldots,\lambda_n\right\},  1\leq i, j\leq n \right\}.
\end{equation*}
%
That is,  the eigenvalues of $(B^\top\otimes A) \mathcal{P}$ is
related to the square roots of the eigenvalues of $\sigma(A^\top
B)$, but from~\eqref{kroeig}, no more information can be used
to decide the positivity or non-negativity of the eigenvalues of
$(B^\top\otimes A) \mathcal{P}$. A question immediately arises as to
whether it is possible to obtain the explicit expression of the
eigenvalues of $(B^\top\otimes A) \mathcal{P}$, provided the
eigenvalues of $A^\top B$ are given. In the following two lemmas,
we first review the periodic QZ decomposition for two matrices and then apply it to discuss the eigenvalues of
$(B^\top\otimes A) \mathcal{P}$.

%
%
%
%


\begin{Lemma}\label{PQZ}
\cite{Bojanczyk1992} Let $A$ and $B$ be two matrices in $\mathbb{R}^{n\times n}$. Then,
 there exist unitary matrices $P,Q\in\mathbb{C}^{n\times n}$ such that
$U_A := PAQ$ and $U_B := Q^H B^\top P^H$   are two upper triangular
matrices.
\end{Lemma}
\begin{Lemma}\label{Lem1}
Let $A$ and $B$ be two matrices in $\mathbb{R}^{m\times n}$. Then
\begin{enumerate}
\item $(B^\top\otimes A)\mathcal{P} =(Q\otimes P^H)(U_A\otimes U_B)\mathcal{P}(Q^H\otimes P)$
\item $\sigma((B^\top\otimes A)\mathcal{P})=
\left \{\lambda_i , \pm\sqrt{\lambda_i\lambda_j} |
 \lambda_i, \lambda_j \in \sigma(A^\top B)  = \left\{ \lambda_1,\ldots,\lambda_n\right\},
1\leq i <  j\leq n \right\}$.
\end{enumerate}
Here,  $\sqrt{z}$ denotes the principal square root of a complex number $z$.
\end{Lemma}

\textbf{Proof}.
Part~1 follows immediately from Lemma~\ref{PQZ} since
$U_A = PAQ$ and $U_B = Q^H B^\top P^H$ for some
unitary matrices $P$ and $Q$, that is,
\begin{eqnarray*}
(B^\top\otimes A)\mathcal{P} &=&(Q\otimes P^H)(U_B\otimes U_A)(P\otimes Q^H)\mathcal{P}\\
&=&(Q\otimes P^H)(U_A\otimes U_B)\mathcal{P}(Q^H\otimes P).
\end{eqnarray*}

Let the diagonal entries of $U_A$ and $U_B$ be denoted by
$\{a_{ii}\}$ and $\{b_{jj}\}$, respectively.
Then, $(U_A \otimes
U_B)$ is an upper triangular matrix with given diagonal entries, specified by $a_{ii}$ and $b_{jj}$.
After multiplying $(U_A \otimes U_B)$ with $\mathcal{P}$ from the
right, the position of the entry $a_{ii}b_{jj}$ is changed to be in the $j+n(i-1)$-th row and the $i+n(j-1)$-th column of the matrix
$(U_A\otimes U_B)\mathcal{P}$. They are then reshuffled by a
sequence of permutation matrices to form a block upper triangular
matrix with diagonal entries arranged in the following order
\begin{eqnarray}\label{BigMat}
&&\left\{ a_{11}b_{11}, \bb
                          0 & a_{11}b_{22} \\
                          a_{22}b_{11} & 0
                      \eb,\ldots,
                      \bb
                          0 & a_{11}b_{nn} \\
                          a_{nn}b_{11} & 0
                      \eb,
                      a_{22}b_{22},
                      \bb
                          0 & a_{22}b_{33} \\
                          a_{33}b_{22} & 0
                      \eb,\right. \nonumber
                      \\
                      &&\left.\ldots,
                      \bb
                          0 & a_{nn}b_{22} \\
                          a_{22}b_{nn} & 0
                      \eb,\ldots,
                      \bb
                          0 & a_{n-1,n-1}b_{nn} \\
                          a_{nn}b_{n-1,n-1} & 0
                      \eb,
                      a_{nn}b_{nn}
 \right\}
\end{eqnarray}
Note that the reshuffling process is not hard to see by following the ordering as used in matrix of size $2$, that is, when $n=2$, $U_A =\bb a_{11} & a_{12} \\ 0 & a_{22} \eb$ and $U_B =\bb b_{11} & b_{12} \\ 0 & b_{22} \eb$, we have
\[ (U_A\otimes
U_B)\mathcal{P} = \bb
a_{11} b_{11}  & a_{12}b_{11}& a_{11} b_{12} &a_{12}b_{12} \\
0 & 0 & a_{11}b_{22} & a_{12}b_{22} \\
0&a_{22}b_{11} &  0  & a_{22} b_{12}\\
0 & 0 & 0 &a_{22} b_{22}
\eb. \]
%
However, it is conceptually simple but operationally tedious to reorder $(U_A\otimes
U_B)\mathcal{P}$ to show this result even for $n=3$ and that will be left as an exercise.

By~\eqref{BigMat}, it can be seen that
\begin{eqnarray*}
\sigma((B^\top\otimes A)\mathcal{P}) &=& \left\{ a_{ii}b_{ii},
\pm\sqrt{a_{ii}a_{jj}b_{ii}b_{jj}}, 1\leq i, j\leq n
\right \} \\
& = &
\left \{\lambda_i , \pm\sqrt{\lambda_i\lambda_j},
1\leq i, j\leq n
\right\}
\end{eqnarray*}
where $\lambda_i = a_{ii}b_{ii} \in \sigma(A^\top B)$  for $1\leq i
\leq n$.

Before demonstrating the unique solvability conditions,
 we need to define that a subset $\Lambda  = \{\lambda_1,\ldots,\lambda_n\}$ of complex numbers is said to be \emph{$\top$-reciprocal free} if and only if whenever $i,j\in\{1,2,\cdots,n\}$, $\lambda_i\neq1/\lambda_j$.
This definition also regards $0$ and $\infty$ as reciprocals of each other. 
 Then, we have the following solvability conditions of~Eq.~\eqref{DTS}.
%
\begin{Theorem}\label{DTSEXIST}
The $\top$-Stein matrix equation~\eqref{DTS} is
uniquely solvable if and only if the following conditions are
satisfied:
\begin{itemize}
\item[a.] The set of $\sigma(A^\top B) \setminus \{-1\}$ is $\top$-reciprocal free.
\item[b.] $-1$ can be an eigenvalue of the matrix $A^\top B$, but must be simple.
\end{itemize}
\end{Theorem}
\textbf{Proof}.
From~\eqref{KronD}, we know that the $\top$-Stein matrix
equation~\eqref{DTS} is  uniquely solvable if and
only if
\begin{equation}\label{SolCond}
1\not\in\sigma((B^\top\otimes A) \mathcal{P}).
\end{equation}
By Lemma~\ref{Lem1},
%
%
if $\lambda \in \sigma(A^\top B)$, then $\dfrac{1}{\lambda}
\not\in\sigma(A^\top B)$.
Otherwise, $1 = \sqrt{\lambda \cdot
\dfrac{1}{\lambda}} \in ((B^\top\otimes A) \mathcal{P})$.
On the other hand, if $-1\in
\sigma(A^\top B)$ and $-1$ is not a simple eigenvalue, then $1\in
\sigma((B^\top\otimes A) \mathcal{P})$. This
verifies~\eqref{SolCond} and the proof of the theorem is complete.
\section{The connection between deflating subspace and Eq.~\eqref{DTS}}
The relationship between solution of matrix equations and the matrix eigenvalue problems has been widely studied in many applications. It is famous that solution of Riccati and polynomial matrix equations can be found by computing invariant subspaces of matrices and deflating subspaces of matrix pencils \cite{Bini2012}. This reality leads us to finding some algorithms for computing solution of Eq.~\eqref{DTS} based on the numerical computation of invariant or deflating subspaces.

Given a pair of $n\times n$ matrices $A$ and $B$, recall that the function $A-\lambda B$ in the variable $\lambda$ is said to be the matrix pencil related to the pair $(A,B)$. For a $k$-dimensional subspace $\mathcal{X}\in\mathbb{C}^{n}$ is called a deflating subspace for the pencil $A-\lambda B$ if there exists a $k$-dimensional subspace $\mathcal{Y}\in\mathbb{C}^n$ such that
\begin{equation*}
A\mathcal{X}\subseteq \mathcal{Y}\mbox{ and }
B\mathcal{X}\subseteq \mathcal{Y},
\end{equation*}
that is,
\begin{equation}\label{eq:invariant}
A {X} = {Y}T_1 \mbox{ and }
B {X} = {Y}T_2,
\end{equation}
where $X, Y \in \mathbb{C}^{n\times k}$ are two full rank matrices whose columns span the spaces $\mathcal{X}$ and $\mathcal{Y}$, respectively, and matrices $T_1, T_2\in\mathbb{C}^{k\times k}$. In particular, if in~\eqref{eq:invariant}, $X=Y$ and $B = T_2 = I$ for an $n\times n$ identity matrix $I$, then we have the simplified formula
\begin{equation}
AX = XT_1.
\end{equation}
Here, the space $\mathcal{X}$ spanned by the columns of the matrix $X$ is called an invariant subspace for $A$, and satisfies
\begin{equation*}
A\mathcal{X} \subseteq \mathcal{X}.
\end{equation*}
One strategy to analyze the eigeninformation is to transform one matrix pencil to its simplified and equivalent form. That is, two matrix pencils $A-\lambda B$ and $\widetilde{A} -\lambda \widetilde{B}$ are said to be equivalent if and only if there exist two nonsingular matrices $P$ and $Q$ such that
 \begin{equation*}
 P(A-\lambda B) Q = \widetilde{A} -\lambda \widetilde{B}.
 \end{equation*}
 In the subsequent discuss, we will use the notion $\sim$ to describe this equivalence relation, i.e., $A-\lambda B \sim \widetilde{A} -\lambda \widetilde{B}$.

Our task in this section is to identify eigenvectors of problem \eqref{eq:invariant}
and then associate these eigenvectors (left and right) with the solution of Eq.~\eqref{DTS}. We begin this analyst
by studying the eigeninformation of two matrices $A$ and $B$, where $A-\lambda B$ is a regular matrix
pencil.

Note that for the ordinary eigenvalue problem,
if the eigenvalues are different then the eigenvectors are linearly independent.
 This property is also true for every regular matrix pencil and is demonstrated as follows.
  For a detailed proof, the reader is referred to~\cite[Theorem 7.3]{Gohberg1982} and~\cite[Theorem 4.2]{Chiang2013}.

\begin{Theorem}\label{lancaster}
Given a pair of $n\times n$ matrix $A$ and $B$, if the matrix pencil $A-\lambda B$ is regular, then its Jordan chains corresponding to all finite and infinite eigenvalues carry the full spectral information about the matrix pencil
and consists of $n$ linearly independent vectors.
 \end{Theorem}
\begin{Lemma}\label{LemmAB}
Let $A-\lambda B\in\mathbb{C}^{n\times n}$ be a regular matrix pencil. Assume that matrices $X_i,Y_i\in\mathbb{C}^{n\times n_i}$, $i=1,2$, are full rank and satisfies the following equations
\begin{subequations}\label{eq:abxy}
\begin{eqnarray}
A X_i &=& Y_i R_i,\\
B X_i &=& Y_i S_i,
\end{eqnarray}
\end{subequations}
where $R_i$ and $S_i$, $i=1,2$, are square matrices of size $n_i\times n_i$.
Then
\begin{itemize}
\item [i)] $R_i-\lambda S_i\in\mathbb{C}^{n_i\times n_i}$ are regular matrix pencils for $i=1,2$.

\item [ii)] if $\sigma(R_1-\lambda S_1)\cap \sigma(R_2-\lambda S_2)=\phi$, then the matrix
$\bb X_1 & X_2\eb\in\mathbb{R}^{n\times (n_1+n_2)}$  is full rank.
\end{itemize}
\end{Lemma}

We also need the following useful lemma.

\begin{Lemma}\label{Chiang Matt}
Given two  regular matrix pencils $A_i-\lambda B_i\in\mathbb{C}^{n_i\times n_i}$, $1\leq i \leq 2$. Consider the following equations with respect to $U,V\in\mathbb{C}^{n_1\times n_2}$
\begin{subequations}
\begin{align}\label{CM}
A_1 U &= V A_2,\\
B_1 U &= V B_2.
\end{align}
\end{subequations}
Then, if $\sigma(A_1-\lambda B_1)\cap\sigma(A_2-\lambda B_2)=\phi$, the equation~\eqref{CM} has the unique solution $U=V=0$.
\end{Lemma}
\textbf{Proof}.
For $n_2=1$, we get

\begin{align*}
A_1 u &= a_2v ,\\
B_1 u &= b_2v ,
\end{align*}

where $a_2,b_2\in\mathbb{C},\,u,v\in\mathbb{C}^{n_1\times 1}$. We may without loss of generality assume that $b_2\neq 0$, then $A_1 u =\frac{a_2}{b_2} B_1 u$ and thus $u=v=0$.
Now, for any $n_2>1$, consider the generalized Schur decomposition of $A_2-\lambda B_2$. We can assume that $A_2=[a_{ij}]$ and $B_2=[b_{ij}]$ are upper triangular matrices (i.e., $a_{ij}=b_{ij}=0,\,1\leq j<i\leq n_2$). Denote that the $i$-th columns of $U$ and $V$ are $u_i$ and $v_i$, respectively.
Thus,
\begin{subequations}\label{CM1}
\begin{align}
A_1 u_i &=\sum\limits_{k=1}^i a_{ki} v_k,\\
B_1 u_i &=\sum\limits_{k=1}^i b_{ki} v_k,
\end{align}
\end{subequations}
for $i=1,2,\ldots n_2$.

If $i=1$, we obtained $u_1=v_1=0$ form the above discussion. Given a integer $i$ such that $1\leq i<n_2$ and assume that $u_k=v_k=0$ for $1\leq k \leq i$. We claim $u_{i+1}=v_{i+1}=0$,
indeed, form \eqref{CM1}, we have
\begin{align*}
A_1 u_{i+1} &= a_{i+1,i+1}v_{i+1} ,\\
B_1 u_{i+1} &= b_{i+1,i+1}v_{i+1}. 
\end{align*}
Again, the result is immediately following the special case $n_2=1$. By mathematical induction we prove this lemma.

\begin{Corollary}
Given $A\in\mathbb{C}^{n\times n}$ and $\Lambda\in\mathbb{C}^{k\times k}$, if $\sigma(A)\cap\sigma(\Lambda)=\phi$. Then the equation with respect to $U\in\mathbb{C}^{n\times k}$
\begin{align*}
AU=U \Lambda
\end{align*}
have the unique solution $U=0$.
\end{Corollary}
Now we have enough tools to analyze the solution of Eq.~\eqref{DTS} associate with some deflating spaces.
 We first establish a important matrix pencil,
let the matrix pencil $\mathcal{M}-\lambda \mathcal{L}$ be defined as
\begin{align}\label{ML2}
\mathcal{M}-\lambda \mathcal{L}:=\bb BA^\top & 0\\-CA^\top & I_n \eb -\lambda \bb I_n & 0 \\ AC^\top & AB^\top \eb \in \mathbb{R}^{2n\times 2n},
\end{align}
it is clear that
\begin{align*}
\sigma(\mathcal{M}-\lambda \mathcal{L})=\sigma(BA^\top) \cup \sigma(I_n-\lambda AB^\top),
\end{align*}
a direct calculation shows that $X$ is a solution of the Eq.~\eqref{DTS} if and only if
\begin{align*}
\mathcal{M} \bb I_n \\ XA^\top \eb &= \bb I_n \\ AX^\top \eb BA^\top,\\
\mathcal{L} \bb I_n \\ XA^\top \eb &= \bb I_n \\ AX^\top \eb
\end{align*}
or if and only if its dual form
\begin{align*}
\bb -AX^\top & I_n \eb\mathcal{M}&=\bb -X A^\top & I_n \eb,\\
 \bb -AX^\top & I_n \eb\mathcal{L}  &= AB^\top \bb -X A^\top & I_n \eb.
\end{align*}

Armed with the property given in Theorem~\ref{lancaster} and Lemma~\ref{Chiang Matt}, we can now attack the problem
of determine how the deflating subspace is related to the solution of Eq.~\eqref{DTS}.
\begin {Theorem}\label{thm2}
Let $A$, $B$ and $C\in\mathbb{R}^{n\times n}$  are given in Eq.~\eqref{DTS},
let us write
\begin{subequations}\label{ML}
 \begin{align}
 \label{sub1}
\mathcal{M}\bb U_1 \\ V_1 \eb=\bb U_2\\V_2 \eb T_1,\\
 \label{sub2}
\mathcal{L}\bb U_1 \\ V_1 \eb=\bb U_2\\V_2 \eb T_2
\end{align}
\end{subequations}
where $\bb U_i\\V_i \eb$ is full rank, $i=1,2$. Assume that the set of $\sigma(BA^\top)$ is $\top$-reciprocal free. Then, we have
\begin{enumerate}
\item $U_1=U_2=0$ if $\sigma (T_1-\lambda T_2)=\sigma(I_n-\lambda AB^\top)$.
\item $U_1$ and $U_2$ are nonsingular  if $T_1-\lambda T_2\sim BA^\top-\lambda I_n$.
 Moreover, if $A$ is nonsingular, then  $X=V_1U_1^{-1}A^{-\top}=U_2^{-\top}V_2^\top A^{-\top}$ is the unique solution  of~Eq.~\eqref{DTS}.
\end{enumerate}
\end {Theorem}
\textbf{Proof}.
From \eqref{ML} we get
\begin{subequations}
\begin{align}
\label{1}
BA^\top U_1 &=U_2 T_1,\\
\label{2}
-C A^\top U_1+V_1 &= V_2 T_1,\\
\label{3}
U_1 &= U_2 T_2,\\
\label{4}
A C^\top U_1+AB^\top V_1 &= V_2 T_2,
\end{align}
\end{subequations}

\begin{itemize}
\item[i)]
It follows from \eqref{1} and \eqref{3} that
since $\sigma(BA^\top-\lambda I_n)\cap\sigma(T_1-\lambda T_2)=\phi$, we have $U_1 = U_2 = 0$ by Lemma~\ref{Chiang Matt}.
\item[ii)]
It can be seen that there exist two nonsingular matrices $U$ and $V$ such that
\begin{subequations}
\begin{align*}
\mathcal{M}\bb 0 \\ U \eb &=\bb 0 \\ V\eb T_2,\\
\mathcal{L}\bb 0 \\ U \eb &=\bb 0 \\ V\eb T_1.  
\end{align*}
\end{subequations}

Hence, together with \eqref{ML} we have
\begin{subequations}
\begin{align*}
\mathcal{M}\bb 0 & U_1\\ U& V_1 \eb &=\bb 0 & U_2\\ V & V_2 \eb
\bb T_2 & 0\\ 0 &  T_1\eb,\\
\mathcal{L}\bb 0 & U_1\\ U& V_1 \eb &=\bb 0 & U_2\\ V & V_2 \eb \bb T_1 & 0\\ 0 &  T_2\eb.
\end{align*}
\end{subequations}

Since $\sigma(\mathcal{M} - \lambda \mathcal{L}) = \sigma(BA^\top - \lambda I_n) \cup \sigma(I_n - \lambda AB^\top)$ and $\sigma(BA^\top - \lambda I_n) \cap \sigma(I_n - \lambda AB^\top)=\phi$, by Theorem~\ref{lancaster} and Lemma~\ref{LemmAB}, the matrix
$\bb 0 & U_1\\ U& V_1 \eb$ is nonsingular. Together with~\eqref{3}, $U_1$ and $U_2$ are nonsingular.

Let $X_i=V_i U_i^{-1}$, $i=1,2$, then form \eqref{2} and \eqref{4}
\begin{align*}
AC^\top +AB^\top X_1 &= V_2 T_2 U_1^{-1}=X_2,\\
-CA^\top +X_1 &=V_2 T_1 U_1^{-1}=X_2 B A^\top,
\end{align*}
or
\begin{align*}
AC^\top +AB^\top X_1 &= X_2,\\
AC^\top +AB^\top X_2^\top &= X_1^\top.
\end{align*}
Since the set of $\sigma(AB^\top)=\sigma(BA^\top)$ is $\top$-reciprocal free, together with
\begin{align*}
X_1^\top-X_2 - AB^\top (X_1^\top-X_2)^\top=0,
\end{align*}
we get $X_1=X_2^\top$. If $A$ is nonsingular, it is easy verify that two matrices $X_1 A^{-\top}$ and $X_2^\top A^{-\top}$
are both satisfying $\top$-Stein equation~Eq.~\eqref{DTS}. The proof of part~(ii) is complete.
\end{itemize}


\begin{Remark}
\begin{itemize}
\item[1.] It is easily seen that $\bb I_n \\ XA^\top \eb$ and $\bb U_1 \\ V_2 \eb$ both span the
unique deflating subspace of $\mathcal{M}-\lambda\mathcal{L}$ corresponding to the set of $\sigma(BA^\top)$. Otherwise, in part~(ii) we know that $T_2$ is nonsingular. We then be able to transform the formulae defined in \eqref{ML} into the generalized eigenvalue problem as follows.
\begin{align*}
\mathcal{M}\bb U_1 \\ V_1\eb= \mathcal{L} \bb U_1 \\ V_1\eb BA^\top.
\end{align*}
That is, some numerical methods for the computation of the eigenspace of $\mathcal{M}-\lambda \mathcal{L}$ corresponding to the set of $\sigma(BA^\top)$ can be designed and solved Eq.~\eqref{DTS}.
\item[2.]
Since the transport of the unique solution $X$ of Eq.~\eqref{DTS} is equal to the unique solution $Y$ of the following matrix equation
\begin{align}\label{dualDTS}
Y=B^\top Y A^\top +C^\top.
\end{align}
Analogous to the consequences of Theorem~\ref{thm2}, The similar results can be obtained with respect to Eq.~\eqref{dualDTS} if $B$ is nonsingular. However, we point out that Eq.~\eqref{DTS} can be solved by computing deflating subspaces of another matrix pencils. For instance we let
\begin{align*}
\mathcal{M}_1-\lambda \mathcal{L}_1:=\bb A^\top B & 0\\-C-AC ^\top B& I_n \eb -\lambda \bb I_n & 0 \\ 0 & AB^\top \eb.
\end{align*}
Assume that the set of $\sigma(BA^\top)$ is $\top$-reciprocal free, it can be shown that $\mathcal{M}_1 \bb I_n \\ X \eb= \mathcal{L}_1 \bb I_n \\ X \eb A^\top B$ and it has similar results as the conclusion of Theorem~\ref{thm2}. The unique solution $X$ of \eqref{DTS} can be found by computing deflating subspaces of the matrix pencil $\mathcal{M}_1-\lambda \mathcal{L}_1$ without the assumption of the singularity of $A$ and $B$.

\end{itemize}
\end{Remark}

\section{Computational methods for solving~Eq.~\eqref{DTS}}
Numerical methods for solving~Eq.~\eqref{DTS} has received great attention in theory and in practice and can be found in~\cite{Wang2007, Su2010} for Krylov subspace methods and in~\cite{Smith68,Penzl99,Zhou2009} for Smith-type iterative methods. In particular, Smith-type iterative methods are only workable in the case $\rho(AB^\top) < 1$, where $\rho(AB^\top)$
denotes the spectral radius of $AB^\top$. In the recent years, a structure algorithm has been studied for Eq.~\eqref{DTS}\cite{Chiang2012} via PQZ decomposition, which consists of transforming  and into Schur form by a PQZ decomposition, and then solving the resulting triangular system by way of back-substitution. In this section, we revisit these numerical methods and point out the advantages and drawbacks of all algorithms.
\subsection{Krylov subspace methods}
Since the $\top$-Stein equation is essentially a linear
system~\eqref{KronD}, we certainly can use Krylov subspace methods to solve~\eqref{KronD}. See, e.g.,
\cite{Wang2007, Su2010}, 
 and the reference cited
therein. The general idea for applying Krylov subspace methods is by
defining the $\star$-Stein operator $\mathcal{T}$ as $\mathcal{T} : X \rightarrow
X - AX^\top B$ and  its adjoint liner operator $\mathcal{T}$ as
$\mathcal{T}^* : Y \rightarrow Y - BY^\top A$ such that $<\mathcal{T}(X), Y> =
<X,\mathcal{T}^*(Y)>$. Here, $X$, $Y\in\mathbb{R}^{m\times n}$ and
the notion $<\cdot,\cdot>$ is denoted as the
Frobenius inner product.
Then, the iterative method based on Krylov subspaces for~Eq.~\eqref{DTS} is as follows.
%

 \begin{itemize}

\item{\bf The conjugate gradient (CG) method}~\cite{Su2010}:
\begin{eqnarray*}
X_{k+1} &=& X_k + \frac{\|R_k\|^2}{\|P_k\|^2}P_k,\\
R_{k+1} &=& C - \mathcal{T}(X_{k+1}) = R_{k} -  \frac{\|R_k\|^2}{\|P_k\|^2}\mathcal{T}(P_k),\\
D_{k+1} &=& \mathcal{T}^*(R_{k+1}) +  \frac{\|R_{k+1}\|^2}{\|R_k\|^2}
D_k,
\end{eqnarray*}
with an initial matrix $X_0$ and the corresponding initial conditions
\begin{equation*}
R_0 = C- \mathcal{T}(X_{0}), \quad
D_0 = \mathcal{T}^*(R_{0}).
\end{equation*}

\end{itemize}

Note that when the solvability conditions of Theorem~\ref{DTSEXIST} are met,
the CG method is guaranteed to converge in a finite number of iterations
for any initial matrix $X_0$.

\subsection{The Bartels-Stewart-like Algorithm~\cite{Bartels1972}}
In this subsection we focus on the discussion of the Bartels-Stewart algorithm, which is known to be a numerical stable algorithm, to solve $\top$-Stein equations. This method is to solve~Eq.~\eqref{DTS} by means of the PQZ decomposition~\cite{Bartels1972}. Its approach has been discussed in~\cite{Chiang2012} and can be summarized as follows. From Lemma~\ref{Lem1}, we know that
there exist two unitary matrices $P$ and $Q$
(see \cite{Bojanczyk1992} for the computation procedure)
 such that
\begin{equation}\label{Ttrial}
PX\overline{Q}  - PAQ\cdot Q^H X^\top P^\top \cdot \overline{P}B \overline{Q}  = PC\overline{Q}
\end{equation}
 With $\widehat{A} = PAQ$ and $\widehat{B}^\top = Q^H B^\top P^H$ being upper-triangular,
the transformed equation looks like
\[
\bb \widehat{X}_{11} & \hat{x}_{12} \\ \hat{x}_{21} &
\hat{x}_{22} \eb  - \bb \widehat{A}_{11} & \hat{a}_{12} \\ 0 & \hat{a}_{22} \eb \bb \widehat{X}_{11}^\top & \hat{x}_{21}^\top \\
\hat{x}_{12}^\top & \hat{x}_{22}^\top \eb \bb \widehat{B}_{11} &  0  \\ \hat{b}_{21} &
\hat{b}_{22} \eb = \bb \widehat{C}_{11} & \hat{c}_{12} \\ \hat{c}_{21} & \hat{c}_{22} \eb
\]
with $\widehat{X} = \bb \widehat{X}_{11} & \hat{x}_{12} \\ \hat{x}_{21} &
\hat{x}_{22} \eb$.
We then have
\begin{eqnarray}
\hat{x}_{22} - \hat{a}_{22} \hat{x}_{22}^\top \hat{b}_{22} &=& \hat{c}_{22}, \label{22x} \\
\hat{x}_{21} - \hat{a}_{22} \hat{x}_{12}^\top \widehat{B}_{11}  &=& \hat{c}_{21} + \hat{a}_{22} \hat{x}_{22}^\top \hat{b}_{21}, \label{21x} \\
\hat{x}_{12} - \widehat{A}_{11} \hat{x}_{21}^\top \hat{b}_{22}  &=& \hat{c}_{12} + \hat{a}_{12} \hat{x}_{22}^\top \widehat{b}_{22}, \label{12x}  ,\\
\widehat{X}_{11} - \widehat{A}_{11} \widehat{X}_{11}^\top \widehat{B}_{11} &=& \widehat{C}_{11} + \widehat{a}_{12}\widehat{x}_{12}^\top \widehat{B}_{11} + \widehat{A}_{11} \widehat{x}_{21}^\top \widehat{b}_{21} + \widehat{a}_{12} \widehat{x}_{22}^\top
\widehat{b}_{21}. \label{11x}
\end{eqnarray}
Thus, the Bartels-Stewart algorithm
can
easily be constructed by first solving $\widehat{x}_{22}$ from~\eqref{22x}, using
$\widehat{x}_{22}$ to obtain $\widehat{x}_{12}$ and $\widehat{x}_{21}$ from~\eqref{21x} and~\eqref{12x}, and then repeating the same discussion  as~\eqref{22x}--\eqref{12x} by
taking advantage of the property of $\widehat{A}_{11}$ and $\widehat{B}_{11}$ being lower triangular matrices from~\eqref{11x}.
\subsection{Smith-type iterative methods}

Originally, Smith-type iterative methods are developed to solve the standard Stein equation
\begin{equation*}
X = \mathcal{A}X\mathcal{B} +\mathcal{C},\quad \mathcal{A},\mathcal{B},\mathcal{C}\in\mathbb{R}^{n\times n}.
\end{equation*}
As mention before, the unknown $X$ is highly related to the generalized eigenspace problems
\begin{subequations}
\begin{align}\label{ML1}
\bb \mathcal{B} & 0\\-\mathcal{C} & I \eb \bb I \\ X \eb  = \bb I & 0 \\ 0 & \mathcal{A} \eb \bb I \\ X \eb \mathcal{B}.
\end{align}
or
\begin{align}\label{MLdual}
\mathcal{A}\bb X & I \eb \bb \mathcal{B} & 0\\ 0 & I\eb=\bb X & I \eb \bb I & 0\\ -\mathcal{C} & \mathcal{A}\eb.
\end{align}
\end{subequations}

Pre-multiplying \eqref{ML1} by the matrix $\bb \mathcal{B} & 0 \\ -\mathcal{A}\mathcal{C} & I \eb$ and post-multiplying \eqref{MLdual} by the matrix $\bb I & 0 \\ -\mathcal{C}\mathcal{B} & \mathcal{A} \eb$
, we get
\begin{align*}
\bb \mathcal{B}^2 & 0\\-\mathcal{C}-\mathcal{A}\mathcal{C}\mathcal{B} & I_n \eb\bb I_n \\ X\eb &=\bb I_n & 0 \\ 0 & \mathcal{A}^2 \eb\bb I_n \\ X\eb \mathcal{B}^2,\\
\mathcal{A}^2 \bb X & I_n\eb \bb \mathcal{B}^2 & 0\\0 & I_n \eb &=\bb X & I_n\eb \bb I_n & 0 \\ -\mathcal{C}-\mathcal{A}\mathcal{C}\mathcal{B} & \mathcal{A}^2 \eb.
\end{align*}
Then, for any positive integer $k>0$, we obtain
\begin{align*}
\bb \mathcal{B}^{2^{k-1}} & 0\\-C_k & I_n \eb\bb I_n \\ X\eb &=\bb I_n & 0 \\ 0 & \mathcal{A}^{2^{k-1}} \eb\bb I_n \\ X\eb \mathcal{B}^{2^{k-1}},\\
\mathcal{A}^{2^{k-1}} \bb X & I_n\eb \bb \mathcal{B}^{2^{k-1}} & 0\\0 & I_n \eb &=\bb X & I_n\eb \bb I_n & 0 \\ -C_k & \mathcal{A}^{2^{k-1}} \eb,
\end{align*}
where the sequence $\{C_k\}$ is defined by
\begin{subequations}\label{Smith}
\begin{align}
C_k&=C_{k-1}+\mathcal{A}^{2^{k-1}}C_{k-1}\mathcal{B}^{2^{k-1}},\quad k\geq 1,\\
C_0&=\mathcal{C}.
\end{align}
\end{subequations}
The explicit expression of $C_k$ is given as following
\begin{align*}
C_k=\sum\limits_{i=1}^{2^k-1} \mathcal{A}^i \mathcal{C} \mathcal{B}^i.
\end{align*}
Under the condition $\rho(\mathcal{A})\rho(\mathcal{B})<1$, it is easy to see that $\{C_k\}$ is convergence, and
\begin{align*}
\limsup\limits_{k\rightarrow\infty}\sqrt[2^k]{\| X-C_k\|} \leq \rho(\mathcal{A})\rho(\mathcal{B}),
\end{align*}
that is, $C_k$ converges quadratically to $X$ as $k\rightarrow\infty$. This iterative method~\eqref{Smith} is called Smith iteration \cite{Smith68}. In recent years, some modified iterative methods are so-called Smith-type iteration, which are based on Smith iteration and improve its speed of convergence. See, e.g., \cite{Zhou2009} and the references cited therein.

Since the condition $\rho(\mathcal{A})\rho(\mathcal{B})<1$ implies that the assumptions of Theorem~\eqref{DTSEXIST} hold, Eq.~\eqref{DTS} is equivalent to Eq.~\eqref{Stein}. We can apply Smith iteration to the Eq.~\eqref{DTS} with the substitution $(\mathcal{A},\mathcal{B},\mathcal{C})=(AB^\top,A^\top B,C+AC^\top B)$.
One possible drawback of the Smith-type iterative
methods is that it cannot always handle the case when there exist
eigenvalues $\lambda, \mu\in\sigma(A^\top B)$ such that
{\bf $\lambda  \mu = -1$} even the unique solution $X$ exist.
Based on the solvable conditions given in this work, it is
possible to develop a specific technique working on the particular
case and it is a subject currently under investigation.

\section{Error analysis}\label{sec:error_ana}
Error analysis is a way for testing the stability of an numerical algorithm and evaluating the accuracy of an approximated solution.
In the subsequent discussion, we want to consider the backward error and perturbation bounds for solving~Eq.~\eqref{DTS}.

As indicated in~\eqref{Ttrial}, matrices $\widehat{A}$ and $\widehat{B}^\top$ are both upper-triangular. We can then apply the error analysis for triangular linear systems in~\cite[Section~3.1]{Golub96}\cite{Higham2002} to obtain

\begin{equation*}
\| \widehat{C} - (\widehat{X} - \widehat{A}\widehat{X}^\top \widehat{B})\|_F
\leq c_{m,n} \mathbf{u}(1+ \|\widehat{A}\|_F\|\widehat{B}\|_F)\|\widehat{X}\|_F,
\end{equation*}
where $c_{m,n}$ is a content depending on the dimensions $m$ and $n$, $\mathbf{u}$ is the unit roundoff.
Since the PQZ decomposition is a stable process, it is true that
\begin{equation}\label{res1}
\| C - ({X} - {A}{X}^\top {B})\|_F
\leq c_{m,n}' \mathbf{u}(1+ \|{A}\|_F\|{B}\|_F)\|{X}\|_F.
\end{equation}
with a modest multiple $c_{m,n}' $.

Note that the inequality of the form~\eqref{res1} can be served as a stopping criterion for terminating iterations generated from Krylov subspace methods~\cite{Wang2007, Su2010} and Smith-type iterative methods~\cite{Smith68,Penzl99,Zhou2009}. In what follows, we shall derive the error associated with numerical algorithms, following the development in~\cite{GhavimiLaub1995, Higham2002}.
\subsection{Backward error}
Like the discussion of ordinary Sylvester equations~\cite{Higham2002}, the normwise backward error of an approximate solution $Y$ of~Eq.~\eqref{DTS} is defined by
\[
\eta(Y) \equiv \min \left\{ \epsilon: Y = (C + \delta C)+(A+\delta A)Y^\top (B + \delta B), \right. \]
\begin{equation}
\left. \|\delta A\|_F \leq \epsilon \alpha, \|\delta B\|_F \leq \epsilon \beta, \|\delta C\|_F \leq \epsilon \gamma \right\},
\end{equation}
where $\alpha \equiv \|A\|_F$, $\beta \equiv \|B\|_F$ and $\gamma \equiv \|C\|_F$.
Let $\mathcal{R} \equiv \delta C + \delta A Y^\top B+ A Y^\top \delta B + \delta A Y^\top \delta B$, which implies that $\mathcal{R} = Y-AY^\top B - C$. It can be seen that the residual $\mathcal{R}$ satisfies
\begin{equation}\label{residual}
\|\mathcal{R}\|_F \leq \eta{(Y)} (\gamma  + \|Y\|_F \alpha\beta (2+\eta(Y))).
\end{equation}
From~\eqref{residual}, we know that a small backward error indeed implies a small relative residual $\mathcal{R}$.
Since the coefficient matrices in~Eq.~\eqref{DTS} include nonlinearity, it appears to be an open problem to obtain the theoretical backward error
with respect to the residual. Again, similar to the Sylvester equation discussed in~\cite[Section~16.2]{Higham2002}, the conditions under which a $\top$-Stein equation has a well-conditioned solution remain unknown.

\subsection{Perturbation bounds}
Consider the perturbed equation
\begin{equation}\label{perturb}
X+\delta X =  (A+\delta A) (X+\delta X)^\top (B+\delta B)+ (C+\delta C).
\end{equation}
Let $S(X) = X - A X^\top B$ be the corresponding $\top$-Stein operator.  We then have
 $S(\delta X) = \delta C + A (X+\delta X)^\top \delta B +
 \delta A (X+\delta X)^\top (B+\delta B)$. With the
 application of norm, it follows that
 \[ \|\delta X\|_F \leq \|S^{-1}\|_F \left\{ \|\delta C\|_F +
\|\delta S\|_F
 (\|X\|_F + \|\delta X\|_F) \right\}, \]
where $\|\delta S\|_F \equiv \|A\|_F \|\delta B\|_F + \|\delta A\|_F (\|B\|_F + \|\delta B\|_F)$.
When $\|\delta S\|_F$ is small enough so that $1 \geq \|S^{-1}\|_F \cdot \|\delta S\|_F$, we can rearrange the above result to
\[ \frac{\|\delta X\|_F}{\|X\|_F} \leq \frac{\|S^{-1}\|_F}{1 - \|S^{-1}\|_F \cdot \|\delta S\|_F} \left( \frac{\|\delta C\|_F}{\|X\|_F} + \|\delta S\|_F \right) \ . \]
With $\|C\|_F = \|S(X)\|_F \leq \|S\|_F \cdot \|X\|_F$ and the condition number $\kappa(S) \equiv \|S\|_F \cdot \|S^{-1}\|_F$, we arrive at the standard perturbation result
\[ \frac{\|\delta X\|_F}{\|X\|_F} \leq \frac{\kappa(S)}{1 - \kappa(S) \cdot \|\delta S\|_F/\|S\|_F} \left( \frac{\|\delta C\|_F}{\|C\|_F} + \frac{\|\delta S\|_F}{\|S\|_F} \right) \ . \]
Thus the relative error in $X$ is controlled by those in $A$, $B$ and $C$, magnified by the condition number $\kappa(S)$.

On the other hand, we can also drop the high order terms in the perturbation to obtain

\begin{equation*}
\delta X -A\delta X^\top B = AX^\top \delta B + \delta A X^\top B+\delta C   .
\end{equation*}
We then rewrite the system in terms of
\begin{equation*}
\mathcal{Q} \mbox{vec}(\delta X) = \left[\!\begin{array}{ccc}(X^\top B)^\top\otimes I_{m} \!&\! I_{n}\otimes (AX^\top) \!&\!  I_{m n} \end{array}\!\right]
\left[\begin{array}{c}\mbox{vec}(\delta A) \\ \mbox{vec}(\delta B) \\ \mbox{vec}(\delta C)\end{array}\right],
\end{equation*}
where $\mathcal{Q} = I_{m n}  - (B^\top\otimes A) \mathcal{P}$. Let $\zeta = \max \left\{\frac{\|\delta A\|_F}{\|A\|_F},
\frac{\|\delta B\|_F}{\|B\|_F}, \frac{\|\delta C\|_F}{\|C\|_F}\right\}$. It can be shown that
\begin{equation}\label{bound}
\frac{\|\delta X\|_F}{\|X\|_F} \leq \sqrt{3}\Psi \zeta,
\end{equation}
where $\Psi = \|\mathcal{Q}^{-1}
\left[\!\begin{array}{ccc}\alpha(X^\top B)^\top\otimes I_{m} \!&\!\beta I_{n}\otimes (AX^\top) \!&\! \gamma I_{m n} \end{array}\!\right]
 \|_2/\|X\|_F.$

A possible disadvantage of the perturbation bound~\eqref{bound}, which ignores the consideration of the underlying structure of the problem, is to overestimate the effect of the perturbation on the data. But this ``universal" perturbation bound is accessible to any given matrices $A$, $B$ and $C$ of~Eq.~\eqref{DTS}.

Unlike the perturbation bound~\eqref{bound}, it is desirable to
obtain a posteriori error bound by assuming
$\delta A = \delta B = 0$ and $\delta C = \widehat{X}-A\widehat{X}^\top B - C$ in~\eqref{perturb}. This assumption gives rise to
\begin{equation}\label{bound3}
\frac{\|\delta X\|_F}{\|X\|_F} \leq
\frac{\|P^{-1}\|_2\|R\|_F}{\|X\|_F}.
\end{equation}
It is true that while doing numerical computation, this bound given in~\eqref{bound3} provides a simpler way for estimating the error
of the solution of~Eq.~\eqref{DTS}.

%
\section{Conclusion}
In this note, we propose a novel approach to the necessary
and sufficient conditions for the unique solvability of the solution $X$ of the
$\top$-Stein equation for square coefficient matrices in terms of the analysis of the spectra $\sigma(A^\top B)$.
Solvability conditions have been derived and algorithms have been proposed in~\cite{Chiang2012, Ikramov2012} by using PQZ decomposition.
On the other hand, one common procedure to solve the Stein-type equations is by means of the invariant subspace method. We believe that our discussion is the first which implements the techniques of the deflating subspace
for solving $\top$-Stein matrix equation and might also gives rise to the possibility of developing an
advanced and effective solver in the future. Also, we obtain the theoretical residual analysis,
 backward error analysis, and perturbation bounds for measuring accurately the error in the computed solution of~Eq.~\eqref{DTS}.


%
%

\section*{Acknowledgement}
The author wish to thank Professor Eric King-wah Chu (Monash University) for many interesting and valuable suggestions on the manuscript. This research work is partially supported by the National Science Council and
the National Center for Theoretical Sciences in Taiwan.

\end{document}